\documentstyle{amsppt}

\define\ep {\varepsilon}
\define\End {\operatorname{End}}
\define\bmo {\operatorname{bmo}}
\define\Lip {\operatorname{Lip}}

\define\zbar {\overline{z}}
\define\zetabar {\overline{\zeta}}
\define\Hbar {\overline{H}}
\define\Ubar {\overline{U}}
\define\Kbar {\overline{K}}

\define\pa {\partial}
\define\RR {\Bbb R}
\define\CC {\Bbb C}
\define\ZZ {\Bbb Z}

\magnification=1200
\vsize=8.5truein

\NoBlackBoxes

\topmatter
\title Integrability of Rough Almost Complex Structures \endtitle

\author C.~Denson Hill and Michael Taylor 
\endauthor

\thanks{The second author was supported by NSF grant DMS-9877077} \endthanks

\endkeywords

 \endsubjclass

\abstract{ We extend the Newlander-Nirenberg theorem to manifolds with
almost complex structures that have somewhat less than Lipschitz regularity.
We also discuss the regularity of local holomorphic coordinates in the
integrable case, with particular attention to Lipschitz almost complex
structures.}
\endabstract

\endtopmatter

\document

$$\text{}$$
{\bf 1. Introduction}
\newline {}\newline

If $\Omega$ is a smooth manifold of dimension $2n$, an almost complex 
structure on $\Omega$ is a section of $\End(T\Omega)$ such that $J^2=-I$.
It is integrable provided there is a local coordinate chart about any
$p\in\Omega$ consisting of holomorphic functions, where $f:\Omega\rightarrow
\CC$ is said to be holomorphic on $\Cal{O}\subset\Omega$ provided
$$
(X+iJX)f=0
\tag{1.1}
$$
on $\Cal{O}$ for all smooth real vector fields $X$ on $\Cal{O}$.  The 
formal integrability condition is that the Lie bracket of complex vector
fields of the form (1.1) continues to have such a form.  Equivalently, one 
forms the Nijenhuis tensor:
$$
N(X,Y)=[X,Y]-[JX,JY]+J[X,JY]+J[JX,Y],
\tag{1.2}
$$
which is linear over $C^1(\Omega)$ and hence defines a tensor field.  The
formal integrability condition is that $N$ vanish.  The Newlander-Nirenberg
theorem [NN] asserts that formal integrability implies integrability.  
In [NN] is was assumed that $J$ had a high degree of smoothness.  Proofs
in [NW] and [Mal] obtained the result for $J\in C^{1+\ep},\ \ep>0$.
See also [Web].
Our goal here is to lower the needed regularity of $J$.  We will establish 
the following.

\proclaim{Theorem 1.1} If $J$ is an almost complex structure on $\Omega$,
satisfying
$$
J\in C^r\cap H^{s,p},\quad r+s>1,\quad sp>2n,
\tag{1.3}
$$
with $0<r<s,\ 1<p<\infty$, then formal integrability implies integrability.
\endproclaim

Here $C^r$ is a H{\"o}lder space and $H^{s,p}$ is an $L^p$-Sobolev space.
For convenience we will work under the hypothesis that also $s<1$.
(But see some comments in \S{5}.)

Our demonstration follows the method of [Mal], but it 
makes use of paraproduct analysis and some elliptic 
regularity results that were perhaps not familiar when the 
papers cited above were written.  As a warm-up, we show that $N$
is well defined when $J$ satisfies (1.3).  Since the components of $N$ 
involve products of components of $J$ and first-order derivatives of
such components, it suffices to show the following.  

\proclaim{Lemma 1.2} Assume $0<r<s<1$ and
$$
u,v\in C^r\cap H^{s,p},\quad r+s>1,\ 1<p<\infty.
\tag{1.4}
$$
Then $u(\pa_j v)$ is a well defined distribution.  Furthermore,
if $u$ and $v$ are defined on a space of dimension $m$,
$$
sp\ge m\Longrightarrow u(\pa_jv)\in C^{r-1}_*\cap H^{s-1,p}.
\tag{1.5}
$$
\endproclaim
\demo{Proof} We use the paraproduct of J.-M.~Bony [Bon] to write
$$
u(\pa_jv)=T_u(\pa_jv)+T_{\pa_jv}u+R(u,\pa_jv).
\tag{1.6}
$$
We have
$$
u\in L^\infty\Longrightarrow T_u\in OP\Cal{B}S^0_{1,1},
\tag{1.7}
$$
the class of operators in $OPS^0_{1,1}$ defined in [Mey].  See also [T],
Chapter 3, whose notation we follow here.
Compare (1.6) with (3.5.1)--(3.5.2) of [T].  Now
$$
T_u\in OP\Cal{B}S^0_{1,1},\ v\in C^r\cap H^{s,p}\Longrightarrow
T_u(\pa_jv)\in C^{r-1}_*\cap H^{s-1,p},
\tag{1.8}
$$
where $C^s_*$ denotes the scale of Zygmund spaces, defined for $s\in\RR$
and coinciding with H{\"o}lder space $C^s$ for $s\in\RR^+\setminus \ZZ^+$.
Next, if $r\in (0,1)$,
$$
\aligned
v\in C^r&\Longrightarrow T_{\pa_jv}\in OP\Cal{B}S^{1-r}_{1,1} \\
&\Longrightarrow T_{\pa_jv}u\in C^{2r-1}_*\cap H^{r+s-1,p},
\endaligned
\tag{1.9}
$$
when $u\in C^r\cap H^{s,p}$; cf.~[T], (3.5.7).  Finally,
$$
v\in C^r\Longrightarrow R_{\pa_jv}\in OPS^{1-r}_{1,1},
\tag{1.10}
$$
where we set $R_g u=R(u,g)$ (cf.~[T], (3.5.11)), and we have
$$
u\in H^{s,p}\Longrightarrow R_{\pa_jv}u\in H^{r+s-1,p},\ {}\
\text{provided }\ r+s-1>0.
\tag{1.11}
$$
For (1.5), it suffices to have $H^{r+s-1,p}\subset C^{r-1}_*$, which holds
if $sp\ge m$.
This proves the lemma.
\enddemo

$\text{}$\newline
{\smc Remark 1.1}. We say more about the operator results used in the proof 
of Lemma 1.2.  A pseudodifferential operator $p(x,D)$ is said to belong to
$OPS^m_{1,\delta}$ provided its symbol $p(x,\xi)$ satisfies
$$
|D_x^\beta D_\xi^\alpha p(x,\xi)|\le C_{\alpha\beta} 
(1+|\xi|)^{m-|\alpha|+\delta|\beta|}.
\tag{1.12}
$$
Operator results on $p(x,D)$ include
$$
p(x,D):H^{r+m,p}\rightarrow H^{r,p},\quad
p(x,D):C^{r+m}_*\rightarrow C^r_*,
\tag{1.13}
$$
assuming $1<p<\infty$.  If $0\le\delta<1$, then (1.13) holds for all
$r\in\RR$.  If $\delta=1$, one needs $r>0$ for (1.13) to hold.  One says
$p(x,D)\in OP\Cal{B}S^m_{1,1}$ provided (1.12) holds with $\delta=1$
and also the partial Fourier transform $\hat{p}(\eta,\xi)$ satisfies
$$
\text{supp}\, \hat{p}(\eta,\xi)\subset \{(\eta,\xi):|\eta|\le a|\xi|\},
\tag{1.14}
$$
for some $a<1$.  For such operators, (1.13) holds for all $r$.

$\text{}$ \newline
{\smc Remark 1.2}.  Another consequence of the operator results mentioned in
Remark 1.1 is that
$$
F\in C^\infty,\ {}\ u\in C^r_*\cap H^{s,p}\Longrightarrow F(u)\in C^r_*
\cap H^{s,p},
\tag{1.15}
$$
given $r,s>0,\ p\in (1,\infty)$.  In fact, tools developed in [Bon] yield
$$
F(u)=F(0)+M(x,D)u,\quad u\in L^\infty\Rightarrow M(x,D)\in OPS^0_{1,1},
\tag{1.16}
$$
and then (1.15) follows from (1.13).  One can have $u=(u_1,\dots,u_K)$,
so in particular $C^r_*\cap H^{s,p}$ contains $uv$ if it contains $u$ and $v$,
and it also contains $u^{-1}$ if this inverse is continuous.

$\text{}$\newline
{\smc Remark 1.3}. Our notation above does not specify domains on which 
elements of $C^r,\ H^{s,p}$, etc., are to be defined.  For example, we could 
take $C^r_*(\RR^m),\ H^{r,p}(\RR^m)$, etc., in (1.13).  Of course a nonzero
constant (perhaps $F(0)$) does not belong to $H^{s,p}(\RR^m)$, but in this 
context we intend to use $H^{s,p}(U)$ and $C^r(\Ubar)$ for some bounded
$U\subset \RR^m$.  The result (1.15) holds in this context.  To justify it
one extends $u$ to an element of $C^r_*\cap H^{s,p}$ on a neighborhood of
$\Ubar$ and makes use of (1.13) in such a situation, as is standard.  
Lemma 1.2 also holds for $u,v\in C^r(\Ubar)\cap H^{s,p}(U)$, and this will be 
useful in \S{3}.

$\text{}$\newline
{\smc Remark 1.4}. Another consequence of Lemma 1.2,
together with Remark 1.2, is that the quantity (1.2) is well
defined for all vector fields $X$ and $Y$ that are regular of class
$C^r\cap H^{s,p}$, under the hypothesis (1.3), and if $N(X,Y)$ vanishes for
all smooth vector fields then it vanishes for all such vector fields.
This class of vector fields is invariant under diffeomorphisms of class
$C^{1+r}\cap H^{1+s,p}$, and this fact will be of use in \S{3}.

$\text{}$\newline
{\smc Remark 1.5}. We dwell a little on the conditions under which the 
regularity hypothesis (1.3) holds.  Note that one special case is
$$
J\in C^{1/2}\cap H^{1/2+\ep,4n},\quad \ep>0.
\tag{1.17}
$$
Now if
$$
J\in C^r,\quad r>\frac{1}{2},
\tag{1.18}
$$
then $J\in H^{1/2+\ep,p}$ for all $p<\infty$ as long as $1/2+\ep<r$,
so Theorem 1.1 applies whenever (1.18) holds.
\newline $\text{}$

The following outlines the rest of this paper.  In \S{2} we rephrase the
problem as an overdetermined system of PDEs on an open set in $\CC^n$,
and give an integrability condition equivalent to $N\equiv 0$.  In \S{3}
we discuss Malgrange's factorization technique, which is to write the 
local coordinate chart $F=(f_1,\dots,f_n)$ as $G\circ H$; in this section 
we construct $H$.  In \S{4} we show that the construction of $G$ is a
consequence of the classical real analytic theory.  All this is in direct 
parallel to [Mal], and the new material in these sections consists of
demonstrations that the various steps work under our weakened regularity 
hypotheses.  In \S{5} we make some concluding comments about the degree 
of regularity of the holomorphic coordinates when hypothesis (1.3) holds,
and we also make special note of the situation when $J$ is Lipschitz.
In such a case, we show the components of $F$ have two derivatives in
$\bmo$.

$$\text{}$$
{\bf 2. Preliminaries}
\newline {}\newline

As in [Mal], we identify a neighborhood of $p$ in $\Omega$ with a
neighborhood of the origin $0$ in $\CC^n$ and arrange that $J(0)$
coincide with the standard complex structure on $\CC^n$.  The task
of solving (1.1) on a neighborhood of $p$ for a family $f_1,\dots,f_n:
\Cal{O}\rightarrow \CC$ forming a local coordinate system becomes that
of solving an overdetermined system of the form
$$
\frac{\pa f_\ell}{\pa \zbar_j}=\sum\limits_k a_{jk}\, 
\frac{\pa f_\ell}{\pa z_k},\quad 1\le j,\ell\le n.
\tag{2.1}
$$
The hypothesis (1.3) is equivalent to
$a_{jk}\in C^r\cap H^{s,p}$,
with $r,s,p$ as in the statement of Theorem 1.1.  
It will be convenient to allow a little wriggle room, when we work with 
(3.9), so we will actually assume
$$
a_{jk}\in C^{r_1}\cap H^{s,p},\quad r_1>r,\quad a_{jk}(0)=0.
\tag{2.2}
$$
This does not affect the validity of Theorem 1.1 as stated.  The condition
$a_{jk}(0)=0$ just reflects our normalization of $J(0)$.

In (2.1), $z_j=x_j+iy_j$ form
the standard coordinates on $\CC^n$ and as usual we set $\pa/\pa\zbar_j=
(1/2)(\pa/\pa x_j+i\pa/\pa y_j)$, etc.  A convenient shorthand is to set
$\pa/\pa \zbar=(\pa/\pa\zbar_1,\dots,\pa/\pa\zbar_n)^t$ (a column vector),
$A_j=(a_{j1},\dots,a_{jn})$ (a row vector), $A=(a_{jk})$, and $F=(f_1,\dots,
f_n)$ (a row vector).  Then (2.1) is written
$$
\frac{\pa F}{\pa \zbar}=A\, \frac{\pa F}{\pa z}.
\tag{2.3}
$$
The integrability condition $N\equiv 0$ becomes
$$
\frac{\pa A_j}{\pa\zbar_k}+A_j\frac{\pa A_k}{\pa z}=
\frac{\pa A_k}{\pa\zbar_j}+A_k\frac{\pa A_j}{\pa z},\quad 1\le j,k\le n.
\tag{2.4}
$$
Note that Lemma 1.2 applies directly to (2.4).  The goal will be to construct 
a solution $F\in C^{1+r}\cap H^{1+s,p}$ to (2.3), with $F$ close to the 
identity map in $C^{1}$-norm, under the assumption that (2.4) holds.
This will provide the desired holomorphic coordinates.

$$\text{}$$
{\bf 3. Malgrange factorization}
\newline {}\newline

A key idea in [Mal] was to produce $F$ in the form
$$
F=G\circ H,
\tag{3.1}
$$
and to apply separate techniques to construct the diffeomorphisms $G$ and $H$.
These factors will be arranged to satisfy
$$
\frac{\pa G}{\pa\zetabar}=B\, \frac{\pa G}{\pa\zeta},
\tag{3.2}
$$
with $\zeta=H(z)$, and
$$
\frac{\pa H}{\pa\zbar}+\frac{\pa\Hbar}{\pa\zbar}(B\circ H)=
A\Bigl[\frac{\pa H}{\pa z}+\frac{\pa\Hbar}{\pa z}(B\circ H)\Bigr],
\tag{3.3}
$$
and furthermore it is arranged that
$$
\sum\limits_j \frac{\pa B_j}{\pa \zeta_j}=0.
\tag{3.4}
$$

Here we need to verify that the construction works under the regularity
hypothesis (2.2).  First, as computed in [Mal], if (3.1) holds with 
$G,H\in C^1$ and if (3.2) holds, then (2.3) is equivalent to (3.3).  
For the next lemma we take $r,s,p$ as in the statement of Theorem 1.1, 
and we assume (2.2) holds.

\proclaim{Lemma 3.1} Assume $H\in C^{1+r}\cap H^{1+s,p}$ and that $H$ is
sufficiently close to the identity in $C^1$-norm, and use (3.3) to define
$B\circ H\in C^r\cap H^{s,p}$, i.e.,
$$
B\circ H=-\Bigl(\frac{\pa\Hbar}{\pa\zbar}-A\frac{\pa\Hbar}{\pa z}\Bigr)^{-1}
\Bigl(\frac{\pa H}{\pa\zbar}-A\frac{\pa H}{\pa z}\Bigr).
\tag{3.5}
$$
Then also $B\in C^r\cap H^{s,p}$.  If $A$ verifies the formal integrability
condition (2.4), then so does $B$, i.e.,
$$
\frac{\pa B_j}{\pa\zetabar_k}+B_j \frac{\pa B_k}{\pa\zeta}=
\frac{\pa B_k}{\pa\zetabar_j}+B_k \frac{\pa B_j}{\pa\zeta},\quad 
1\le j,k\le n.
\tag{3.6}
$$
\endproclaim
\demo{Proof}  As (2.4) restates the formal integrability of $J$, i.e., the
vanishing of (1.2) for all vector fields $X,Y\in C^{r}\cap H^{s,p}$, the
system (3.6) restates the formal integrability of $\widetilde{J}$, given by
$$
\widetilde{J}(\zeta)=DH(z)\, J(z)\, DH(z)^{-1},\quad \zeta=H(z).
$$
Note that under the stated regularity hypotheses on $J$ and $H$ we have
$\widetilde{J}\in C^r\cap H^{s,p}$.  Now the equivalence of the formal 
integrability of $J$ and $\widetilde{J}$ just amounts to the coordinate
invariance of (1.2).
\enddemo

The next result extends the scope of the lemma on p.~294 of [Mal].
\proclaim{Lemma 3.2} Given $\ep,\delta>0$, one can find a ball $U$ about
$0\in\CC^n$ and
$$
H\in C^{1+r}(\Ubar)\cap H^{1+s,p}(U),
\tag{3.7}
$$
satisfying 
$$
H(0)=0,\quad \|H-id\|_{C^{1+r}(\Ubar)}<\delta,
\tag{3.8}
$$
and such that $B\in C^r(\Ubar)\cap H^{s,p}(U)$, defined by (3.5), satisfies
$|B(0)|<\ep$ and also satisfies (3.4).
\endproclaim

In our proof, following [Mal], we find it convenient to fix $U$ and dilate
$A$, setting $A_t(z)=A(tz)$.  Let $E=B\circ H$, so
$$
E=-\Bigl(\frac{\pa\Hbar}{\pa\zbar}-A_t\frac{\pa\Hbar}{\pa z}\Bigr)^{-1}
\Bigl(\frac{\pa H}{\pa\zbar}-A_t \frac{\pa H}{\pa z}\Bigr).
\tag{3.9}
$$
The map $\Phi(H,t)=E$ has the mapping property
$$
\Phi:\Cal{B}^{r,s,p}(\delta)\times[0,1]\longrightarrow 
C^r(\Ubar)\cap H^{s,p}(U),
\tag{3.10}
$$
where $\Cal{B}^{r,s,p}(\delta)$ consists of $H\in C^{1+r}(\Ubar)\cap 
H^{1+s,p}(U)$ satisfying (3.8).  Then both $\Phi$ and $D_H\Phi$ are 
continuous, where we regard our Banach spaces as real Banach spaces.
An application of the chain rule gives
$$
\frac{\pa B_j}{\pa \zeta_j}\circ H=\Bigl(\frac{\pa K}{\pa\zeta_j}\circ H
\Bigr)\frac{\pa E_j}{\pa z}+\Bigl(\frac{\pa\Kbar}{\pa\zeta_j}\circ H
\Bigr)\frac{\pa E_j}{\pa\zbar},\quad K=H^{-1}.
\tag{3.11}
$$
One can express $(\pa K/\pa\zeta_j)\circ H$ and $(\pa\Kbar/\pa\zeta_j)\circ H$
in terms of the $z$- and $\zbar$-derivatives of $H$ and $\Hbar$, 
using the identity
$$
(DK)\circ H(x)=DH(x)^{-1}
$$
of $(2n)\times (2n)$ real matrices.  It follows 
via Lemma 1.2 (extended to function spaces on bounded domains) that
$$
\Psi(H,t)=\sum\limits_j \frac{\pa B_j}{\pa\zeta_j}\circ H
\tag{3.12}
$$
defines a map
$$
\Psi:\Cal{B}^{r,s,p}(\delta)\times [0,1]\longrightarrow 
C^{r-1}_*(\Ubar)\cap H^{s-1,p}(U),
\tag{3.13}
$$
which is continuous, along with $D_H\Psi$.  Note that $\Psi(id,0)=0$.
A calculation gives
$$
D_H\Psi(id,0)h=\sum\limits_j \frac{\pa^2h}{\pa z_j\pa\zbar_j}=\frac{1}{4}
\Delta h.
\tag{3.14}
$$
This map has a right inverse
$$
\widetilde{\Cal{G}}h=4(\Cal{G}h-\Cal{G}h(0)),
\tag{3.15}
$$
where $\Cal{G}$ denotes the solution operator to
$$
\Delta v=h\ \text{ on }\ U,\quad v\bigr|_{\pa U}=0.
\tag{3.16}
$$
This has the crucial mapping properties
$$
\Cal{G}:C^{r-1}_*(\Ubar)\rightarrow C^{r+1}(\Ubar),\quad
\Cal{G}:H^{s-1,p}(U)\rightarrow H^{s+1,p}(U),
\tag{3.17}
$$
valid for $0<r,s<1,\ 1<p<\infty$.  These properties of $\Cal{G}$ 
can be established by extending $h$ 
to a neighborhood of $\Ubar$, using local regularity, and 
then using well known regularity of harmonic functions on $\Ubar$ with
boundary values in $C^{r+1} (\pa U)$ or $B^{s+1-1/p}_{p,p}(\pa U)$,
respectively.  The implicit function theorem yields for small $t>0$ a 
solution to $\Psi(H,t)=0$ close to the solution $H_0(z)=z$ to $\Psi(H_0,0)=0$.
This proves the lemma.

$\text{}$ \newline
{\smc Remark 3.1}. It is for the continuity in $t$ in (3.10) and (3.13) that
we need $r_1>r$ in (2.2).  Also for this reason we need $sp>2n$ in (1.3),
rather than the weaker inequality that suffices for (1.5).

$$\text{}$$
{\bf 4. Reduction to the analytic case}
\newline {}\newline

At this point it remains to construct a diffeomorphism $G\in C^{1+r}$
satisfying (3.2), where $B\in C^r\cap H^{s,p}$ satisfies (3.4) and (3.6).
The key observation of [Mal] is that this forces $B$ to be real analytic.
Our final task is to verify that this works under our weaker regularity
hypothesis.

\proclaim{Lemma 4.1} If $r,s,p$ satisfy the hypotheses of Theorem 1.1 and
(3.4) and (3.6) hold for $B\in C^r\cap H^{s,p}$, then $B$ is real analytic,
on a neighborhood of $0\in \CC^n$.
\endproclaim
\demo{Proof} Since $|B(0)|<\ep$,  the system (3.4), (3.6) is an overdetermined
elliptic system for $B$ on a neighborhood of the origin.  Once we show
$B\in C^\infty$, the real analyticity is classical.  The smoothness of $B$
follows from Theorem 2.2.E of [T], but for the reader's convenience we sketch
a proof.

The overdetermined elliptic system (3.4), (3.6) has the form
$$
LB+K(B,\nabla B)=0,
\tag{4.1}
$$
where $L$ is a first-order linear operator (with constant coefficients, in
this case) and $K(B,\nabla B)$ is bilinear in its arguments.  Ellipticity 
near $z=0$ follows from $|B(0)|<\ep$.  Using the sort of symbol smoothing
discussed in Chapter 3 of [T], we can write
$$
K(B,\nabla B)=M^\#(x,D)B+M^b(x,D)B,
\tag{4.2}
$$
and, by (3.3.25) of [T], given $r>0$,
$$
B\in C^r\Longrightarrow M^\#(x,\xi)\in S^1_{1,\delta},\ {}\
M^b(x,\xi)\in S^{1-r\delta}_{1,1}.
\tag{4.3}
$$
Here $\delta\in (0,1)$ is picked in advance.  Also we have $L+M^\#(x,D)\in
OPS^1_{1,\delta}$ elliptic, with left parametrix $E\in OPS^{-1}_{1,\delta}$.
Hence, if $\delta$ is close enough to $1$, we have
$$
B=-E\, M^b(x,D)B,\quad \text{mod}\, C^\infty.
\tag{4.4}
$$
and
$$
\aligned
B\in C^r\cap H^{s,p},\ r+s>1&\Longrightarrow M^b(x,D)B\in H^{s+r\delta-1,p} \\
&\Longrightarrow B\in H^{s+r\delta,p}.
\endaligned
\tag{4.5}
$$
We can iterate (4.4)--(4.5) arbitrarily often to obtain $B\in C^\infty$.
\enddemo

Having Lemma 4.1, the endgame is that given in [Mal].  Since $B$ is real 
analytic, the Cartan-Kahler theorem implies (3.2) is solvable, for a real
analytic diffeomorphism $G$, given the integrability condition (3.6).

$$\text{}$$
{\bf 5. Further regularity results}
\newline {}\newline

Here we make note of the regularity of the map $F$ in terms of the 
hypothesized regularity of the almost complex structure $J$.  Here is one
result.

\proclaim{Proposition 5.1} Under the hypothesis (1.3) on $J$, when $N
\equiv 0$, then
$$
F\in C_*^{r+1}\cap H^{s+1,p}.
\tag{5.1}
$$
\endproclaim

Once we show that $H$ has this regularity, then (5.1) will follow for
$F=G\circ H$, as in Remark 1.2.  This degree of regularity for $H$ was
indicated in Lemma 3.2, but to establish this lemma we raised the 
regularity assumed on $J$, in (2.2), so, shifting back, at this point 
we merely have $H\in C^{\rho+1}_*\cap H^{s+1,p}$, for $r-\rho=\ep>0$,
arbitrarily small.  We will be able to go from here to (5.1) via some
elliptic regularity.  In fact, the formulas (3.9) and (3.11) show that
the condition (3.4) yields for $H$ an elliptic system of the form
$$
\sum\limits_j a_j(\nabla H)\pa_j b_j(\nabla H,A\nabla H)=0,
\tag{5.2}
$$
with $a_j$ and $b_j$ smooth in its arguments and
$$
A\in C^r\cap H^{s,p},\quad H\in C^{\rho+1}\cap H^{s+1,p},
\tag{5.3}
$$
with $r-\rho=\ep>0$, arbitrarily small.  
We continue to impose the conditions in 
(1.3) on $r,s,p$, and for simplicity we continue to assume $0<r<s<1$.
The deduction that $H\in C^{r+1}$ is not a standard elliptic regularity 
result, but we can bring paraproduct techniques to bear to prove it.

To be sure, passing from $F\in C^{r+1-\ep}$ to $F\in C^{r+1}$ is a small 
thing.  In fact, the main point of this section is to establish the
next regularity result.

\proclaim{Proposition 5.2} If $J$ is Lipschitz and $N\equiv 0$, then
$$
\pa^2F\in\bmo.
\tag{5.4}
$$
\endproclaim

Here $\bmo$ denotes the localized John-Nirenberg space.
What we know from the results of \S\S{3--4} is that, under this hypothesis,
$H\in C^{1+r}$ and hence $F\in C^{1+r}$, for all $r<1$.  Again it remains
to establish that $\pa^2H\in\bmo$ when $H$ satisfies an elliptic system 
of the form (5.2) and we know that
$$
A\in\Lip,\quad H\in C^{1+\rho},\quad \forall\ \rho<1,
\tag{5.5}
$$
and (equivalently) $H\in H^{s+1,p},\ \forall\ s<1,p<\infty$.

To begin the proof, the paradifferential calculus described in Remark 1.2
gives
$$
b_j(\nabla H,A\nabla H)=B_{j1}\nabla H+B_{j2}A\nabla H\ \text{ mod }C^\infty,
\tag{5.6}
$$
where we have $B_{j\nu}\in OPS^0_{1,1}$, and furthermore, upon choosing
$\delta\in (0,1)$, we can write
$$
\gathered 
\nabla H,A\nabla H\in C^\rho\Longrightarrow B_{j\nu}=B^\#_{j\nu}+
B^b_{j\nu}, \\ B^\#_{j\nu}\in OPS^0_{1,\delta},\quad B^b_{j\nu}
\in OPS^{-\rho\delta}_{1,1}.
\endgathered
\tag{5.7}
$$
Cf.~[T], Proposition 3.1.D.  Hence we have
$$
B^b_{j1}\nabla H,\ B^b_{j2}A\nabla H\in C^{\rho+\rho\delta}_*\cap
H^{s+\rho\delta,p}.
\tag{5.8}
$$
Thus (5.2) yields
$$
\sum\limits_j a_j(\nabla H)\pa_j(B^\#_{j1}\nabla H+B^\#_{j2}A\nabla H)=f_1,
\tag{5.9}
$$
with
$$
f_1=-\sum a_j (\nabla H)\pa_j(B^b_{j1}\nabla H+B^b_{j2}A\nabla H)
\in C^{\rho+\gamma-1}_*\cap H^{s+\gamma-1,p},
\tag{5.10}
$$
for some $\gamma>0$, whose specific formula we do not need.  Next, we
analyze the product in (5.9) in terms of paraproducts, as in (1.6), 
obtaining
$$
\sum\limits_j T_{a_j(\nabla H)} \pa_j(B^\#_{j1}\nabla H+
B^\#_{j2}A\nabla H)=f_2,
\tag{5.11}
$$
with
$$
\aligned
f_2&=f_1-\sum T_{\psi_j}a_j(\nabla H)-\sum R(a_j(\nabla H),\psi_j), \\
\psi_j&=\pa_j(B^\#_{j1}\nabla H+B^\#_{j2}A\nabla H).
\endaligned
\tag{5.12}
$$
We have $\psi_j\in C^{\rho-1}_*\cap H^{s-1,p}$.  Then, via (1.9),
$$
T_{\psi_j}a_j(\nabla H)\in C^{2\rho-1}_*\cap H^{\rho+s-1,p},
\tag{5.13}
$$
while, as in (1.11), given that $\rho$ is so close to $r$ that $\rho+s>1$,
$$
R(a_j(\nabla H),\psi_j)\in H^{\rho+s-1,p}\subset C^{\rho+\gamma-1}_*,\
\text{ if }\ sp>2n,
\tag{5.14}
$$
where again $\gamma$ is a positive number that we need not compute.
In summary,
$$
f_2\in C^{\rho+\gamma-1}_*.
\tag{5.15}
$$
Next, we have
$$
T_{a_j(\nabla H)}=P_j^\#+P_j^b,\quad
P_j^\#\in OPS^0_{1,\delta},\quad P_j^b\in OP\Cal{B}S^{-\rho\delta}_{1,1}.
\tag{5.16}
$$
Then (5.11) yields
$$
\sum\limits_j P_j^\# \pa_j(B^\#_{j1}\nabla H+B^\#_{j2}A\nabla H)=f_3,
\tag{5.17}
$$
where, with $\psi_j$ as in (5.12),
$$
f_3=f_2-\sum P^b_j\psi_j\in C^{\rho+\gamma-1}_*,
$$
and again $\gamma>0$.

For the next step, we have
$$
A\nabla H=T_A\nabla H+T_{\nabla H}A+R(A,\nabla H)=T_A \nabla H+g.
\tag{5.18}
$$
This is where the regularity of $A$ crucially affects the regularity of $H$.
Given $\nabla H\in C^\rho$, $0<\rho<1$, we have, for $r>0$,
$$
A\in C^r_*\Rightarrow g\in C^r_*,\quad A\in\Lip\Rightarrow \nabla g\in\bmo.
\tag{5.19}
$$
Now (5.17) yields
$$
\sum\limits_j P^\#_j\pa_j(B^\#_{j1}\nabla H+B^\#_{j2}T_A\nabla H)=f_4,
\tag{5.20}
$$
where, with $g$ as in (5.18)--(5.19),
$$
f_4=f_3-\sum P^\#_j\pa_jB^\#_{j2}g.
\tag{5.21}
$$
Hence, given $0<\rho<r,\ \rho+\gamma>r$ (which can be assumed in this 
context),
$$
A\in C^r_*\Rightarrow f_4\in C^{r-1}_*,\quad 
A\in\Lip\Rightarrow f_4\in\bmo.
\tag{5.22}
$$
Furthermore, we have
$$
T_A=Q^\#+Q^b,\quad Q^\#\in OPS^0_{1,\delta},\quad Q^b\in OP\Cal{B}
S^{-r\delta}_{1,1},
\tag{5.23}
$$
where we can take $r=1$ if $A\in\Lip$.  Then (5.20) yields
$$
\sum\limits_j P^\#_j\pa_j(B^\#_{j1}\nabla H+B^\#_{j2}Q^\#\nabla H)=f_5,
\tag{5.24}
$$
with
$$
f_5=f_4-\sum P^\#_j\pa_jB^\#_{j2}Q^b \nabla H,
\tag{5.25}
$$
which has as much regularity as $f_4$.  Now the left side of (5.24) is of the 
form $PH$, where $P\in OPS^2_{1,\delta}$ is elliptic, and hence has a 
parametrix $E\in OPS^{-2}_{1,\delta}$.  Hence
$$
H=Ef_5,\ \text{ mod }\ C^\infty,
\tag{5.26}
$$
and, by (5.22) and the analogue for $f_5$, given $0<r<1$,
$$
A\in C^r\Rightarrow H\in C^{r+1},\quad A\in\Lip\Rightarrow \pa^2H\in\bmo.
\tag{5.27}
$$
This completes the proof of Propositions 5.1 and 5.2.

$$\text{}$$
{\bf References}
\newline $\text{}$
\roster
\item"[Bon]" J.-M.~Bony, Calcul symbolique et propagation des singulariti\'es
pour des \'equations aux d\'eriv\'ees nonlin\'eaires, Ann. Sci. Ecole Norm.
Sup. 14 (1981), 209--246.
\item"[Mal]" B.~Malgrange, Sur l'int\'egrabilit\'e des structures 
presque-complexes, Symposia Math., Vol. II (INDAM, Rome, 1968),
Academic Press, London, 1969, pp.~289--296.
\item"[Mey]" Y.~Meyer, Remarques sur un th\'eor\`eme de J.M.~Bony,
Rend. del Circolo mat. di Palermo (suppl. II:1) (1981), 1--20.
\item"[NN]" A.~Newlander and L.~Nirenberg, Complex coordinates in almost 
complex manifolds, Ann. of Math. 65 (1957), 391--404.
\item"[NW]" A.~Nijenhuis and W.B.~Woolf, Some integration problems in
almost-complex manifolds, Ann. of Math. 77 (1963), 424--489.
\item"[T]" M.~Taylor, Pseudodifferential Operators and Nonlinear PDE,
Birkhauser, Boston, 1991.
\item"[Web]" S.~Webster, A new proof of the Newlander-Nirenberg theorem, 
Math. Zeit. 201 (1989), 303--316.

\endroster

$$\text{}$$
Math. Dept., State University of New York, Stony Brook NY  11794
\newline
Math. Dept., University of North Carolina, Chapel Hill NC 27599

\enddocument